\documentclass[11pt,twoside]{article}
\usepackage{a4}
\usepackage{amssymb,amsmath,amsthm,latexsym}
\usepackage{amsfonts}
  \usepackage{amsfonts}
\usepackage{graphicx}

 \def\ab{{\bold a}}
 
 \def\xb{{\bold x}}

 \def\opn#1#2{\def#1{\operatorname{#2}}} % to make operators
 \opn\chara{char} \opn\length{\ell} \opn\pd{pd} \opn\rk{rk}
 \opn\projdim{proj\,dim} \opn\injdim{inj\,dim} \opn\rank{rank}
 \opn\depth{depth} \opn\grade{grade} \opn\height{height}
 \opn\embdim{emb\,dim} \opn\codim{codim}
 
 \opn\Tr{Tr} \opn\bigrank{big\,rank}
 \opn\superheight{superheight}\opn\lcm{lcm}
 \opn\trdeg{tr\,deg}%\emph{
 \opn\reg{reg} \opn\lreg{lreg} \opn\ini{in} \opn\lpd{lpd}
 \opn\size{size} \opn\sdepth{sdepth}
 \opn\link{link}\opn\fdepth{fdepth}\opn\lex{lex}
 \opn\div{div} \opn\Div{Div} \opn\cl{cl} \opn\Cl{Cl}
 \opn\Spec{Spec} \opn\Supp{Supp} \opn\supp{supp} \opn\Sing{Sing}
 \opn\Ass{Ass} \opn\Min{Min}\opn\Mon{Mon}
 \opn\Ann{Ann} \opn\Rad{Rad} \opn\Soc{Soc}
 \opn\Im{Im} \opn\Ker{Ker} \opn\Coker{Coker} \opn\Am{Am}
 \opn\Hom{Hom} \opn\Tor{Tor} \opn\Ext{Ext} \opn\End{End}
 \opn\Aut{Aut} \opn\id{id}
 
 \opn\nat{nat}
 \opn\pff{pf}%   \pf exists already
 \opn\Pf{Pf} \opn\GL{GL} \opn\SL{SL} \opn\mod{mod} \opn\ord{ord}
 \opn\Gin{Gin} \opn\Hilb{Hilb}\opn\sort{sort}
 \opn\Tot{Tot}
 \opn\aff{aff} \opn
 \con{conv} \opn\relint{relint} \opn\st{st}
 \opn\lk{lk} \opn\cn{cn} \opn\core{core} \opn\vol{vol}
 \opn\link{link} \opn\star{star}\opn\lex{lex}\opn\set{set}
 \opn\gr{gr}
 
 \def\pot#1#2{#1[\kern-0.28ex[#2]\kern-0.28ex]}

 \opn\dirlim{\underrightarrow{\lim}}
 \opn\inivlim{\underleftarrow{\lim}}
 \let\union=\cup
 \let\sect=\cap

 \def\Implies{\ifmmode\Longrightarrow \else
         \unskip${}\Longrightarrow{}$\ignorespaces\fi}
 \def\implies{\ifmmode\Rightarrow \else
         \unskip${}\Rightarrow{}$\ignorespaces\fi}
 \def\iff{\ifmmode\Longleftrightarrow \else
         \unskip${}\Longleftrightarrow{}$\ignorespaces\fi}

 \let\:=\colon

\newtheorem{theorem}{Theorem}[section]

\newtheorem{corollary}[theorem] {Corollary}
\newtheorem{definition}[theorem]{Definition}
\newtheorem{example}[theorem]{Example}

\newtheorem{remark}[theorem]{Remark}

\vspace{5cm}

\begin{document}

  \label{'ubf'}
\setcounter{page}{1}                                 %Put here the starting page number

\begin{center}
{
       {\Large \textbf { Unmixedness of generalized Veronese bi-type ideals
    % Put the title of the paper here
                               }
       }
\\

\medskip

{ Monica La Barbiera$^{1}$ and Roya Moghimipor$^{2}$}\\
{\footnotesize $^{1}$Department of Mathematics and Informatics, University of Messina, Viale Ferdinando
Stagno d’Alcontres 31, 98166 Messina, Italy}\\
{\footnotesize  e-mail: {\it monicalb@unime.it}}\\

{\footnotesize $^{2}$Department of Mathematics, Safadasht Branch, Islamic
Azad University, Tehran, Iran}\\

{\footnotesize e-mail: {\it roya\_moghimipour@yahoo.com}}
}
\end{center}

\thispagestyle{empty}

\begin{abstract}
{\footnotesize  In this paper, some algebraic invariants of generalized Veronese bi-type ideals are computed. We characterize the unmixed generalized Veronese bi-type ideals and we give a description of their associated prime ideals.
}
 \end{abstract}



\section{Introduction}
Let $K$ be a field and $K[x_1,\dots,x_n]$ the polynomial
ring in $n$ variables over $K$ with each $x_i$ of degree $1$.
Let $I\subset S$ be a monomial ideal and $G(I)$ its unique
minimal monomial generators.

Let $K$ be a field and  $S = K[x_1, \ldots,  x_n, y_1, \ldots, y_m]$ be the polynomial ring over $K$ in the variables $x_i$ and $y_j$.
In \cite{L3} the first author introduced a class of monomial ideals of $S$, so-called Veronese bi-type ideals. They are an extension
of the ideals of Veronese type (\cite{V}) in a polynomial ring in two sets of variables.
More precisely, the ideals of Veronese bi-type are monomial ideals of $S$ generated in the same degree $q$:
$L_{q,s}=\sum _{k+r=q}I_{k,s}J_{r,s}$, with $k,r\geq1$, $s\leq q$, where
$I_{k,s}$ is the Veronese type ideal generated in degree $k$ by the set
\[
\{x_{1}^{a_{1}}\dots x_{n}^{a_{n}}| \sum_{i=1}^{n}a_{i}=k, 0\leq a_{i} \leq s , s \in \{1,\dots,k\}\}
\]
and $J_{r,s}$ is the Veronese type ideal generated in degree $r$ by the set
\[
\{y_{1}^{b_{1}}\dots y_{m}^{b_{m}}| \sum_{j=1}^{m}b_{j}=r, 0\leq b_{j} \leq s , s \in \{1,\dots,r\}\}.
\]
For $s=2$ the Veronese bi-type ideals are the ideals of the walks of a bipartite graph with loops (\cite{L3}).
The first author \cite{B} studied the combinatorics
of the integral closuer and the normality of $L_{q,2}$. More in general, in \cite{L3} the same problem is studied for $L_{q,s}$ for all $s$.
A great deal of knowledge on the Veronese bi-type ideal is accumulated in several papers \cite{{BPR}, {B}, {L3}, {L2}, {L}}.

Now we consider the polynomial ring $T$ over $K$ in the variables
\[
x_{11},\ldots,x_{1m_1},x_{21},\ldots,x_{2m_2},\ldots,x_{n1},\ldots,x_{nm_n}.
\]

In this paper, we introduce the concept of generalized Veronese bi-type ideals.
The concept generalized Veronese bi-type ideal generalizes the concept of Veronese bi-type ideals.
Let $t, s, q_{1}, \dots ,q_{n}$ be non negative integers with $s\leq t$ and $\sum_{i=1}^{n}q_{i}=t$, $q_{1}, \dots ,q_{n}\geq1$.
The {\em ideals of generalized Veronese bi-type} of degree $t$ are the monomial ideals of $T$
\[
L^{*}_{t,s}=\sum_{s\leq t, \sum_{i=1}^{n}q_{i}=t}L_{1,q_{1},s}\dots L_{n,q_{n},s},
\]
where the ideals $L_{i,q_{i},s}$ are Veronese type ideals of degree $q_{i}$ generated by the monomials $x_{i1}^{a_{i1}} \ldots x_{im_{i}}^{a_{im_{i}}}$
with $\sum_{j=1}^{m_{i}}a_{ij}=q_{i}$ and $0 \leq a_{ij}\leq s$ for $i=1,\dots,n$.
When $s=2$, the generalized Veronese bi-type ideals arise from n-partite graphs with loops, the so-called strong quasi-n-partite graphs.
A graph $G$ with loops is said to be quasi-n-partite if its vertex set $V=V_1\union V_2\union \cdots\union V_n$ and $V_i=\{x_{i1},\ldots,x_{im_i}\}$ for $i=1,\ldots,n$, every edge joins a vertex of $V_{i}$ with a vertex of $V_{i+1}$, and some vertices in $V$ have loops.
A quasi-n-partite graph is called {\em strong} if it is a complete n-partite graph and all its vertices have loops.
A strong quasi-n-partite graph on vertices $x_{11},\ldots,x_{1m_1},\ldots,x_{n1},\ldots,x_{nm_n}$ will be denoted by $\mathcal{K}'_{m_1,\dots,m_n}$.

The present paper is organized as follows. In Section \ref{two} unmixed ideals of generalized Veronese bi-type are classified and the generalized ideals
associated to the walks of special n-partite graphs, described by the generalized Veronese bi-type ideals
\[
L^{*}_{t,2}=\sum_{\sum_{i=1}^{n}q_{i}=t}L_{1,q_{1},2}\dots L_{n,q_{n},2},
\]
are considered in \cite{Vi, Vl}.
Furthermore we investigate some algebraic invariants of $T/I(L^{*}_{t,s})$.

In Section \ref{three} we give in  Theorem \ref{Ass} a description of the associated prime ideals of generalized Veronese bi-type ideals.

In Section \ref{four} the toric ideal $I(L^{*}_{t,s})$ of the monomial subring $K[L^{*}_{t,s}]\subset T$ is studied.
Let $L^{*}_{t,s}=(f_1,\dots,f_p)$ and $K[L^{*}_{t,s}]$ be the K-algebra spanned by $f_1,\dots,f_p$.
There is a graded epimorphism of K-algebras: $\varphi: R=K[t_1,\dots,t_p]\rightarrow K[L^{*}_{t,s}]$ induced by $\varphi(t_i)=f_i$, where $R$ is a polynomial ring graded by $\deg(t_i)=\deg(f_i)$. Let $I(L^{*}_{t,s})$ be the toric ideal of $K[L^{*}_{t,s}]$, that is the kernel of $\varphi$.
In Corollary \ref{Groebner} we show that $I(L^{*}_{t,s})$ has a quadratic Groebner basis and as a consequence the K-algebra $K[L^{*}_{t,s}]$ is Koszul.

\section{Generalized Veronese bi-type ideals}
\label{two}
Let $S=K[x_1,\ldots,x_n]$ be the polynomial ring over a field $K$ in the variables $x_1,\ldots,x_n$, and let $I\subset S$ be a monomial ideal with $I\neq S$ whose minimal set of generators  is $G(I)=\{\xb^{\ab_1},\ldots, \xb^{\ab_r}\}$. Here $\xb^{\ab_i}=x_1^{\ab_i(1)}x_2^{\ab_i(2)}\cdots x_n^{\ab_i(n)}$ for $\ab_i=(\ab_i(1),\ldots,\ab_i(n))\in\mathbb{Z}^{n}_{+}=\{\mathbf{u}=(u_1,\ldots,u_n)\in \mathbb{Z}^{n}: u_i \geq 0 \}$.
We consider the polynomial ring $T$ over $K$ in the variables
\[
x_{11},\ldots,x_{1m_1},x_{21},\ldots,x_{2m_2},\ldots,x_{n1},\ldots,x_{nm_n}.
\]

Now we introduce a class of monomial ideals of $T$, the so-called generalized Veronese bi-type ideals, which are an extension of the ideals of Veronese bi-type
introduced in \cite{B}.
Let $t, s, q_{1}, \dots ,q_{n}$ be non negative integers with $s\leq t$ and $\sum_{i=1}^{n}q_{i}=t$, $q_{1}, \dots ,q_{n}\geq1$.
The {\em ideals of generalized Veronese bi-type} of degree $t$ are the monomial ideals of $T$
\[
L^{*}_{t,s}=\sum_{s\leq t, \sum_{i=1}^{n}q_{i}=t}L_{1,q_{1},s}\dots L_{n,q_{n},s},
\]
where the ideals $L_{i,q_{i},s}$ are Veronese type ideals of degree $q_{i}$ generated by the monomials $x_{i1}^{a_{i1}} \ldots x_{im_{i}}^{a_{im_{i}}}$
with $\sum_{j=1}^{m_{i}}a_{ij}=q_{i}$ and $0 \leq a_{ij}\leq s$ for $i=1,\dots,n$.

\begin{remark}
{\em
In general $L_{i,q_{i},s}\subseteq L_{i,q_{i}}$ for all $i=1,\dots,n$, where $L_{i,q_{i}}$ is the {\em Veronese ideal} of degree $q_{i}$ generated by all the monomials in the variables $x_{i1},\ldots,x_{im_i}$ of degree $q_{i}$ \cite{S, V}.

One has $L_{i,q_{i},s}=L_{i,q_{i}}$ for any $q_{i}\leq s$. If $s=1$, $L_{i,q_{i},1}$ is the {\em squarefree Veronese ideal} of degree $q_{i}$
generated by all the squarefree monomials in the variables $x_{i1},\ldots,x_{im_i}$ of degree $q_{i}$.
}
\end{remark}

\begin{example}
\label{L_{2,2}}
{\em
Let $T=K[x_{11},x_{12},x_{21},x_{22}]$ be a polynomial ring.

\item{(1)}
$L^{*}_{2,2}=L_{1,1,2}L_{2,1,2}=L_{1,1}L_{2,1}=(x_{11}x_{21},x_{11}x_{22},x_{12}x_{21}, x_{12}x_{22})$.
\begin{eqnarray*}
(2) L^{*}_{4,2}
&=&L_{1,3,2}L_{2,1,2}+L_{1,1,2}L_{2,3,2}+L_{1,2,2}L_{2,2,2}=L_{1,3,2}L_{2,1}+L_{1,1}L_{2,3,2}+L_{1,2}L_{2,2}\\
&=&(x_{11}^{2}x_{12}x_{21},x_{11}^{2}x_{12}x_{22}, x_{11}x_{12}^{2}x_{21},x_{11}x_{12}^{2}x_{22}, x_{11}x_{21}^{2}x_{22}, x_{12}x_{21}^{2}x_{22}, x_{11}x_{21}x_{22}^{2},\\
&&x_{12}x_{21}x_{22}^{2}, x_{11}^{2}x_{21}^{2}, x_{11}^{2}x_{21}x_{22}, x_{11}^{2}x_{22}^{2},x_{12}^{2}x_{21}^{2},x_{12}^{2}x_{22}^{2},x_{12}^{2}x_{21}x_{22},
x_{11}x_{12}x_{21}^{2},\\
&&x_{11}x_{12}x_{22}^{2}, x_{11}x_{12}x_{21}x_{22}).
\end{eqnarray*}
}\end{example}

Next we investigate algebraic invariants of $T/L^{*}_{t,s}$.
It would be appropriate to recall the definition of the Castelnuovo-Mumford regularity.
We refer the reader to \cite{E} for further details on the subject.

Let $M$ be a finitely generated graded $S$-module.
The {\em Castelnuovo-Mumford regularity} (or simply the regularity) of $M$ is defined as
\[
\reg(M) := \max_{i,j}\{j - i : \beta_{i,j}(M) \neq 0\},
\]
where $\beta_{i,j}(M)=\dim_{K}(Tor_{i}(K,M))_{j}$ denotes the $ij$-th graded Betti number of $M$.

\begin{theorem}
\label{reg}
Let $L^{*}_{t,s}$ be a generalized Veronese bi-type ideal of $T$. Then
\[
\reg(T/L^{*}_{t,s})=t-1.
\]
\end{theorem}

\begin{proof}
Let $t, s, q_{1}, \dots ,q_{n}$ be non negative integers with $s\leq t$ and $\sum_{i=1}^{n}q_{i}=t$, $q_{1}, \dots ,q_{n}\geq1$.
Let
\[
L^{*}_{t,s}=\sum_{s\leq t, \sum_{i=1}^{n}q_{i}=t}L_{1,q_{1},s}\dots L_{n,q_{n},s}
\]
be a generalized Veronese bi-type ideal, where the ideals $L_{i,q_{i},s}$ are Veronese type ideals of degree $q_{i}$ generated by the monomials $x_{i1}^{a_{i1}} \ldots x_{im_{i}}^{a_{im_{i}}}$ with $\sum_{j=1}^{m_{i}}a_{ij}=q_{i}$ and $0 \leq a_{ij}\leq s$.
Then
\[
\reg(T/L^{*}_{t,s})=\max \{\deg f \mid f \quad \text{minimal generator of} \quad L^{*}_{t,s}\}-1=t-1,
\]
as desired.
\end{proof}

\begin{example}
\label{L_{3,2}}
{\em
Let $T=K[x_{11},x_{12},x_{21},x_{22},x_{31},x_{32}]$ be a polynomial ring.
Let
\begin{eqnarray*}
L^{*}_{3,2}
&=&L_{1,1,2}L_{2,1,2}L_{3,1,2}\\
&=&L_{1,1}L_{2,1}L_{3,1}\\
&=&(x_{11}x_{21}x_{31},x_{11}x_{21}x_{32},x_{11}x_{22}x_{31},x_{11}x_{22}x_{32},x_{12}x_{21}x_{31},x_{12}x_{21}x_{32}, x_{12}x_{22}x_{31},\\
&&x_{12}x_{22}x_{32})
\end{eqnarray*}
be a generalized Veronese bi-type ideal of $T$.
It follows from Theorem \ref{reg} that $\reg( T/L^{*}_{3,2})=3-1=2$.
}
\end{example}

A {\em vertex cover} of $L^{*}_{t,s}$ is a subset $W$ of
\[
\{x_{11},\ldots,x_{1m_1},x_{21},\ldots,x_{2m_2},\ldots,x_{n1},\ldots,x_{nm_n}\}
\]
such that each $u\in G(L^{*}_{t,s})$ is divided by some variables of $W$.
Such a vertex cover is called {\em minimal}
if no proper subset of $W$ is vertex cover.
We denote the minimal cardinality of the vertex covers of $L^{*}_{t,s}$ by $h(L^{*}_{t,s})$.

\begin{theorem}
\label{dim}
Let $L^{*}_{t,s}$ be a generalized Veronese bi-type ideal of $T$. Then one has:
\item[(a)] if $2\leq t< s(m_1+\dots+m_n)-r$ for $r=1,\dots,s-1$, then
\[
\dim(T/L^{*}_{t,s})=m_1+\dots+m_n-\min\{m_1,\dots,m_n\}.
\]

\item[(b)] if $t=s(m_1+\dots+m_n)-r$ for $r=1,\dots,s-1$, then
$\dim(T/L^{*}_{t,s})=m_1+\dots+m_n-1$.
\end{theorem}

\begin{proof}
(a) By the structure of $G(L^{*}_{t,s})$ the minimal vertex covers of $L^{*}_{t,s}$ are $W_{i}=\{x_{i1},\dots,x_{im_{i}}\}$ for $i=1,\dots,n$.
The minimal cardinality of the vertex covers of $L^{*}_{t,s}$ is $h(L^{*}_{t,s})=\min\{m_1,\dots,m_n\}$. Therefore,
\begin{eqnarray*}
\dim(T/L^{*}_{t,s})
&=& m_1+\dots+m_n- h(L^{*}_{t,s})\\
&=&m_1+\dots+m_n-\min\{m_1,\dots,m_n\}.
\end{eqnarray*}

(b) Suppose that $t=s(m_1+\dots+m_n)-r$ for $r=1,\dots, s-1$. Thus the minimal cardinality of the vertex covers is $h(L^{*}_{t,s})=1$, being $W=\{x_{11}\}$ a minimal vertex cover of $L^{*}_{t,s}$ by construction.
\end{proof}

\begin{example}
\label{dimension}
{\em
Let
\begin{eqnarray*}
L^{*}_{3,2}
&=&(x_{11}x_{21}x_{31},x_{11}x_{21}x_{32},x_{11}x_{22}x_{31},x_{11}x_{22}x_{32},x_{12}x_{21}x_{31},x_{12}x_{21}x_{32}, x_{12}x_{22}x_{31},\\
&&x_{12}x_{22}x_{32})
\end{eqnarray*}
be a generalized Veronese bi-type ideal of $T=K[x_{11},x_{12},x_{21},x_{22},x_{31},x_{32}]$.
Hence, by applying Theorem \ref{dim} we obtain that $\dim(T/L^{*}_{3,2})=6-2=4$.
}
\end{example}

A monomial ideal is said to be {\em unmixed}
if all its minimal vertex covers have the same cardinality.
We recall the one-to-one correspondence between the minimal
vertex covers of an ideal and its minimal prime ideals.
Thus $P$ is a minimal prime ideal of $L$ if and only if
$P=(\mathcal{A})$ for some minimal vertex cover $\mathcal{A}$ of $L$.

In the following, we classify the unmixed ideals of generalized Veronese bi-type.

\begin{theorem}
\label{unmixed1}
Let $L^{*}_{t,s}$ be a generalized Veronese bi-type ideal of $T$ with $2\leq t< s(m_1+\dots+m_n)-r$ for $r=1,\dots,s-1$. Then
$L^{*}_{t,s}$ is unmixed if and only if $m_{1}=m_{2}=\dots=m_{n}$.
\end{theorem}

\begin{proof}
Let $T=K[x_{11},\ldots,x_{1m_1},x_{21},\ldots,x_{2m_2},\ldots,x_{n1},\ldots,x_{nm_n}]$. Suppose that $2\leq t< s(m_1+\dots+m_n)-r$ for $r=1,\dots,s-1$.
It then follows from Theorem ~\ref{dim} that the minimal cardinality of the vertex covers of $L^{*}_{t,s}$ is $h(L^{*}_{t,s})=\min\{m_1,\dots,m_n\}$.
Therefore, all the minimal vertex covers have the same cardinality if and only if $m_{1}=m_{2}=\dots=m_{n}$.
\end{proof}

\begin{example}
\label{unmixedL_{3,2}}
{\em
Let $T=K[x_{11},x_{12},x_{21},x_{22}]$ be a polynomial ring. Let
\begin{eqnarray*}
L^{*}_{3,2}&=&(x_{11}^{2}x_{21},x_{11}^{2}x_{22},x_{11}x_{12}x_{21}, x_{11}x_{12}x_{22}, x_{12}^{2}x_{21}, x_{12}^{2}x_{22},x_{11}x_{21}^{2}, x_{11}x_{21}x_{22},x_{11}x_{22}^{2},\\
&& x_{12}x_{21}^{2},x_{12}x_{21}x_{22},x_{12}x_{22}^{2}),
\end{eqnarray*}
be a generalized Veronese bi-type ideal of $T$.
The minimal vertex covers are:
$W_{1}=\{x_{11},x_{12}\}$; $W_{2}=\{x_{21},x_{22}\}$.
Therefore, $h(L^{*}_{3,2})=|W_{1}|=|W_{2}|=2$, and hence $L^{*}_{3,2}$ is unmixed by Theorem \ref{unmixed1}.
}
\end{example}

\begin{theorem}
\label{unmixed2}
Let $L^{*}_{t,s}$ be a generalized Veronese bi-type ideal of $T$ with $t=s(m_1+\dots+m_n)-r$ for $r=1,\dots,s-1$. Then
$L^{*}_{t,s}$ is unmixed.
\end{theorem}

\begin{proof}
For all $i=1,\ldots,n$ and $j=1,\ldots,m_i$, one has $W_{i}=\{x_{ij}\}$ are the minimal vertex covers of $L^{*}_{t,s}$ by construction.
\end{proof}

\begin{example}
\label{unmixedL_{11,3}}
{\em
Let
\begin{eqnarray*}
L^{*}_{11,3}=(x_{11}^{3}x_{12}^{3}x_{21}^{3}x_{22}^{2},x_{11}^{3}x_{12}^{3}x_{21}^{2}x_{22}^{3},x_{11}^{3}x_{12}^{2}x_{21}^{3}x_{22}^{3},
x_{11}^{2}x_{12}^{3}x_{21}^{3}x_{22}^{3}),
\end{eqnarray*}
be a generalized Veronese bi-type ideal of $T=K[x_{11},x_{12},x_{21},x_{22}]$.
The minimal vertex covers of $L^{*}_{11,3}$ are:

$W_{1}=\{x_{11}\}$; $W_{2}=\{x_{12}\}$; $W_{3}=\{x_{21}\}$; $W_{4}=\{x_{22}\}$. Therefore, $h(L^{*}_{11,3})=|W_{i}|=1$ for all $i=1,2,3,4$, and hence $L^{*}_{11,3}$ is unmixed.
}
\end{example}

As an application, we consider ideals arising from graph theory.

A graph $G$ consists of a finite set $V=\{x_1,\dots,x_n\}$ of vertices and a collection $E(G)$ of subsets of $V$, that consists of pairs $\{x_i,x_j\}$, for some
$x_i,x_j\in V$.

A graph $G$ has loops if it is not requiring $x_i\neq x_j$ for all edges $\{x_i,x_j\}$ of $G$. Then the edge $\{x_i,x_i\}$ is said a loop of $G$.

A graph $G$ with loops is called {\em complete} if each pair $\{x_i,x_j\}$ is an edge of $G$ for all $x_i,x_j\in V.$

We observe that the ideals of generalized Veronese bi-type can be associated to graphs with loops.

\begin{definition}
{\em
A graph $G$ with loops is said to be {\em quasi-n-partite} if its vertex set $V=V_1\union V_2\union \cdots\union V_n$ and $V_i=\{x_{i1},\ldots,x_{im_i}\}$ for $i=1,\ldots,n$, every edge joins a vertex of $V_{i}$ with a vertex of $V_{i+1}$, and some vertices in $V$ have loops.
}
\end{definition}

\begin{definition}
{\em
A quasi-n-partite graph $G$ is called {\em strong} if it is a complete n-partite graph and all its vertices have loops.
}
\end{definition}

A strong quasi-n-partite graph on vertices $x_{11},\ldots,x_{1m_1},\ldots,x_{n1},\ldots,x_{nm_n}$ will be denoted by $\mathcal{K}'_{m_1,\dots,m_n}$.

Let $G$ be a graph with loops in each of its $n$ vertices. A {\em walk} of {\em length} $t$ in $G$ is an alternating sequence
\[
w=\{v_{i_{0}},l_{i_{1}},v_{i_{1}},l_{i_{2}},\dots, v_{i_{t-1}}, l_{i_{t}},v_{i_{t}}\},
\]
where $v_{i_{0}}$ or $v_{i_{g}}$ is a vertex of $G$ and
$l_{i_{g}}=\{v_{i_{g-1}},v_{i_{g}}\}$, $g=1,\dots,t$, is either the edge joining $v_{i_{g-1}}$ and $v_{i_{g}}$ or a loop if $v_{i_{g-1}}=v_{i_{g}}, 1\leq i_{0}\leq i_{1}\leq \dots\leq i_{t}\leq n$.

\begin{example}
{\em
Let $\mathcal{K}'_{n,m}$ be a strong quasi-bipartite graph on vertices $x_{1},\ldots,x_{n}$, $y_{1},\ldots,y_{m}$.
A walk of length 2 in $\mathcal{K}'_{n,m}$ is
\[
\{x_i,l_{i},x_i,l_{ij},y_{j}\} \quad  \text{or} \quad \{x_i,l_{ij},y_j,l_j,y_j\}
\]
where $l_{i}=\{x_i,x_i\}$, $l_{j}=\{y_j,y_j\}$ are loops, and $l_{ij}$
is the edge joining $x_i$ and $y_j$. Because $\mathcal{K}'_{n,m}$ is bipartite, any walk in it have not the edges $\{x_{i_{h}},x_{i_{k}}\}$, $i_h\neq i_k$, and
$\{y_{j_{h}},y_{j_{k}}\}$, $j_h\neq j_k$.
}
\end{example}

Let $G$ be a graph with loops. The {\em generalized graph ideal} $I_t(G)$ associated to $G$ is the ideal of the polynomial ring $T$ generated by all the monomials of degree $t\geq 3$ corresponding to the walks of length $t-1$. Thus, the variables in each generator of $I_t(G)$ have at most degree $2$.

Now let $\mathcal{K}'_{m_1,\dots,m_n}$ be a strong quasi-n-partite graph with vertex set $V=V_1\union V_2\union \cdots\union V_n$ and $V_i=\{x_{i1},\ldots,x_{im_i}\}$ for $i=1,\ldots,n$. For this graph we have
\[
I_{t}(\mathcal{K}'_{m_1,\dots,m_n})=L^{*}_{t,2}=\sum_{\sum_{i=1}^{n}q_{i}=t}L_{1,q_{1},2}\dots L_{n,q_{n},2},
\]
for $t\geq3$.

\begin{remark}
{\em
If $t=2$, the ideal $L^{*}_{t,2}$ does not describe the edge ideal
\[
I(\mathcal{K}'_{m_1,\dots,m_n})=I_{2}(\mathcal{K}'_{m_1,\dots,m_n})
 \]
of a strong quasi-n-partite graph. Let $\mathcal{K}'_{2,2}$ be the strong quasi-bipartite graph on vertices $x_{11},x_{12},x_{21},x_{22}$,
then
\[
I(\mathcal{K}'_{2,2})=(x_{11}x_{21},x_{11}x_{22},x_{12}x_{21},x_{12}x_{22},x_{11}^{2},x_{12}^{2},x_{21}^{2},x_{22}^{2}),
\]
but
$L^{*}_{2,2}=(x_{11}x_{21},x_{11}x_{22},x_{12}x_{21},x_{12}x_{22})$. Therefore, $I(\mathcal{K}'_{2,2})\neq L^{*}_{2,2}$.
}
\end{remark}

\begin{example}
{\em
Let $T=K[x_{11},x_{12},x_{21},x_{22}]$ be a polynomial ring and $\mathcal{K}'_{2,2}$ be the strong quasi-bipartite graph on vertices $x_{11},x_{12},x_{21},x_{22}$. Then
\begin{eqnarray*}
I_{4}(\mathcal{K}'_{2,2})
&=&L_{1,3,2}L_{2,1}+L_{1,1}L_{2,3,2}+L_{1,2}L_{2,2}\\
&=&(x_{11}^{2}x_{12}x_{21},x_{11}^{2}x_{12}x_{22}, x_{11}x_{12}^{2}x_{21},x_{11}x_{12}^{2}x_{22}, x_{11}x_{21}^{2}x_{22}, x_{12}x_{21}^{2}x_{22}, x_{11}x_{21}x_{22}^{2},\\
&&x_{12}x_{21}x_{22}^{2}, x_{11}^{2}x_{21}^{2}, x_{11}^{2}x_{21}x_{22}, x_{11}^{2}x_{22}^{2},x_{12}^{2}x_{21}^{2},x_{12}^{2}x_{22}^{2},x_{12}^{2}x_{21}x_{22},
x_{11}x_{12}x_{21}^{2},\\
&&x_{11}x_{12}x_{22}^{2}, x_{11}x_{12}x_{21}x_{22}).
\end{eqnarray*}
}
\end{example}

The following result classifies the ideals $I_{t}(G)$ that are unmixed.

\begin{theorem}
\label{unmixedgraph}
Let $T=K[x_{11},\ldots,x_{1m_1},x_{21},\ldots,x_{2m_2},\ldots,x_{n1},\ldots,x_{nm_n}]$.
\item[(a)] If $2\leq t<2(m_{1}+\dots+m_{n})-1$, then $L^{*}_{t,2}$ is unmixed if and only if $m_1=\dots=m_n$.

\item[(b)] If $t=2(m_{1}+\dots+m_{n})-1$, then $L^{*}_{t,2}$ is unmixed.
\end{theorem}

\begin{proof}
The assertion follows by Theorem ~\ref{unmixed1} and ~\ref{unmixed2}.
\end{proof}

\begin{example}
\label{k2,2,2}
{\em
Let $T=K[x_{11},x_{12},x_{21},x_{22},x_{31},x_{32}]$ be a polynomial ring and $\mathcal{K}'_{2,2,2}$ be the strong quasi-3-partite graph on vertices $x_{11},x_{12},x_{21},x_{22},x_{31},x_{32}$. Therefore,
\begin{eqnarray*}
I_{3}(\mathcal{K}'_{2,2,2})
&=&(x_{11}x_{21}x_{31},x_{11}x_{21}x_{32},x_{11}x_{22}x_{31},x_{11}x_{22}x_{32},x_{12}x_{21}x_{31},x_{12}x_{21}x_{32}, x_{12}x_{22}x_{31},\\
&&x_{12}x_{22}x_{32}).
\end{eqnarray*}
Then Theorem \ref{unmixedgraph} implies that $I_{3}(\mathcal{K}'_{2,2,2})$ is unmixed.
}
\end{example}

\section{Associated prime ideals of generalized Veronese bi-type ideals}
\label{three}
In this section we want to determine the associated prime ideals of generalized Veronese bi-type ideals.
Let $S=K[x_{1},\dots,x_{n}]$ be the polynomial ring over a field $K$ in the variables $x_1,\ldots,x_n$ with the maximal ideal $\mathfrak{m}=(x_1,\dots,x_n)$, and let $I\subset S$ be a monomial ideal with $I\neq S$ whose minimal set of generators  is $G(I)=\{\xb^{\ab_1},\ldots, \xb^{\ab_m}\}$.

A prime ideal $P \subseteq S$ is an {\em associated prime} of $I$ if there exists an element $a\in S$ such that $I : (a)=P$.
The set of associated primes of an ideal $I$ in a ring $S$ is to be denoted by $\Ass_{S}(S/I)$.
Next we consider the polynomial ring $T$ over $K$ in the variables
\[
x_{11},\ldots,x_{1m_1},x_{21},\ldots,x_{2m_2},\ldots,x_{n1},\ldots,x_{nm_n}.
\]
Let $\mathcal{F}\subseteq \{1,2,\ldots,m_{1}+\dots+m_{n}\}$, where $m_{1}+\dots+m_{n}$ is the number of the variables of the polynomial ring $T$.
For a subset $\mathcal{F}$ we denote by $\mathcal{P}_{\mathcal{F}}$ the prime ideal of $T$ generated by the variables whose index is in $\mathcal{F}$.

\begin{theorem}
\label{Ass}
Let $L^{*}_{t,s}$ be a generalized Veronese bi-type ideal of $T$.
\[
\mathcal{P}_{\mathcal{F}}\in \Ass_{T}(T/L^{*}_{t,s})\Longleftrightarrow |\mathcal{F}|\leq r+1,
\]
for $r=s(m_{1}+\dots+m_{n})-t$, $r=1,\dots, s-1$.
\end{theorem}

\begin{proof}
We replace the set of variables $\{x_{11},\dots, x_{1m_{1}}\}$ with $\{y_{1},\dots,y_{m_{1}}\}$ and $\{x_{21},\dots$, $x_{2m_{2}}\}$ with
$\{y_{m_{1}+1},\dots, y_{m_{1}+m_{2}}\}$ and so on up to $\{x_{n1},\dots,x_{nm_{n}}\}$ with
\[
\{y_{m_{1}+\dots +m_{n-1}+1},\dots,y_{m_{1}+\dots+m_{n}}\}.
\]

Suppose that $\mathcal{P}_{\mathcal{F}}\in \Ass_{T}(T/L^{*}_{t,s})$. Thus there exists a monomial $f\notin L^{*}_{t,s}$ such that $L^{*}_{t,s}:f=\mathcal{P}_{\mathcal{F}}$. We show that we can choose such a monomial $f$ of degree $t-1$ such that $L^{*}_{t,s}:f=\mathcal{P}_{\mathcal{F}}$.

Assume that $f\notin L^{*}_{t,s}$, $L^{*}_{t,s}:f=\mathcal{P}_{\mathcal{F}}$, $\deg(f)\geq t$ and
$f=y_{1}^{b_{1}}\dots y_{m_{1}+\dots+m_{n}}^{b_{m_{1}+\dots+m_{n}}}$. Then there exists $z\in \{1,2,\ldots,m_{1}+\dots+m_{n}\}$ such that $b_{z}>s$.
Since $L^{*}_{t,s}:f=\mathcal{P}_{\mathcal{F}}$, it follows that $b_{r}f\in L^{*}_{t,s}$ for all $r\in \mathcal{F}$ and $b_{r}f\notin L^{*}_{t,s}$ for all $r\notin \mathcal{F}$. Furthermore, for all $r\in \mathcal{F}$ there exists a monomial $u_{r}\in G(L^{*}_{t,s})$ such that $u_{r}|(y_{r}f)$.
Being $f\notin L^{*}_{t,s}$, this fact means that, for all $r\in \mathcal{F}$, the variable $y_{r}$ appears in $u_{r}$ with exponent $b_{r}+1$.
Then $b_{r}<s$ for all $r\in \mathcal{F}$, and hence $z\notin \mathcal{F}$.

Now we claim that: 1) $\overline{f}=f/y_{z}\notin L^{*}_{t,s}$ and 2) $L^{*}_{t,s}:\overline{f}=\mathcal{P}_{\mathcal{F}}$.

The first fact follows from that $f\notin L^{*}_{t,s}$ and $b_{z}-1\geq s$. For the second assertion we proceed as follows.
$L^{*}_{t,s}:\overline{f}\subseteq L^{*}_{t,s}:f$ because $\overline{f}$ divides $f$. Then $L^{*}_{t,s}:\overline{f}\subseteq \mathcal{P}_{\mathcal{F}}$, being
$\mathcal{P}_{\mathcal{F}}=L^{*}_{t,s}:f$. Since $b_{z}-1\geq s$, it follows that $u_{r}$ divides $y_{r}f/y_{z}$ for all $r\in \mathcal{F}$, hence
$y_{r}\in L^{*}_{t,s}: (f/y_{z})$ for all $r\in \mathcal{F}$. Thus $\mathcal{P}_{\mathcal{F}}\subseteq L^{*}_{t,s}: (f/y_{z})$.
It follows the other inclusion $\mathcal{P}_{\mathcal{F}}\subseteq L^{*}_{t,s}:\overline{f}$. Then $\mathcal{P}_{\mathcal{F}}=L^{*}_{t,s}:\overline{f}$.
After a finite number of these reductions, we find $f\notin L^{*}_{t,s}$ of degree $t-1$ such that $\mathcal{P}_{\mathcal{F}}= L^{*}_{t,s}:f$.
It then follows that $fy_{r}\in L^{*}_{t,s}$ for all $r\in \mathcal{F}$ and $fy_{r}\notin L^{*}_{t,s}$ for all $r\notin \mathcal{F}$.
More precisely $b_{r}+1\leq s$ for all $r\in \mathcal{F}$, and $b_{r}\leq s$ for all $r\notin \mathcal{F}$.
Therefore, $b_{r}=s$ for all $r\notin \mathcal{F}$, and hence $f=\prod_{r\in \mathcal{F}}y_{r}^{b_{r}} \prod_{r\notin \mathcal{F}}y_{r}^{s}$ with
$0\leq b_{r}< s$ for all $r\in \mathcal{F}$. We have
\[
\deg(\prod_{r\notin \mathcal{F}}y_{r}^{s})=s(m_{1}+\dots+m_{n}-|\mathcal{F}|)=q.
\]
Thus
\begin{eqnarray*}
s(m_{1}+\dots+m_{n})
&\geq& (\sum_{r\in \mathcal{F}}b_{r}+1)+q\\
&=&\sum_{r\in \mathcal{F}}b_{r}+|\mathcal{F}|+q\\
&=&\sum_{r\in \mathcal{F}}b_{r}+q+|\mathcal{F}|\\
&=&\deg(f)+|\mathcal{F}|\\
&=&t-1+|\mathcal{F}|.
\end{eqnarray*}
Conversely, let $|\mathcal{F}|\leq r+1$, for $r=s(m_{1}+\dots+m_{n})-t$ , $r=1,\dots,s-1$, that is $|\mathcal{F}|\leq s(m_{1}+\dots+m_{n})-t+1$.
Furthermore, in these hypotheses one has $s(m_{1}+\dots+m_{n}-|\mathcal{F}|)\leq t-1$. In fact,
\[
s(m_{1}+\dots+m_{n}-|\mathcal{F}|) \leq s(m_{1}+\dots+m_{n})-r-1;
\]
thus $s|\mathcal{F}|\geq r+1$ that is true for $r=1,\dots,s-1$.
We assume that $t=s(m_{1}+\dots+m_{n})-r$ for $r=1,\dots,s-1$. Then, for any monomial $u\in G(L^{*}_{t,s})$,
there exists an integer $p\in \mathcal{F}$ such that $y_{p}$ divides $u$. Thus $L^{*}_{t,s} \subset \mathcal{P}_{\mathcal{F}}$.
The condition $s(m_{1}+\dots+m_{n})\geq t-1+|\mathcal{F}|$ implies that
\[
(s-1)|\mathcal{F}|+s(m_{1}+\dots+m_{n}-|\mathcal{F}|)\geq t-1,
\]
which together with $s(m_{1}+\dots+m_{n}-|\mathcal{F}|)\leq t-1$ shows that there exists an integer $d_{r}< s$, for all $r\in \mathcal{F}$ such that
\[
d_{r}|\mathcal{F}|+s(m_{1}+\dots+m_{n}-|\mathcal{F}|)=t-1.
\]
Then the monomial $f=\prod_{r\in \mathcal{F}}y_{r}^{d_{r}}\prod_{r\notin \mathcal{F}}y_{r}^{s}$
has degree $t-1$. Thus $f\notin L^{*}_{t,s}$ and as a consequence $\mathcal{P}_{\mathcal{F}}\subseteq L^{*}_{t,s}:f$.

Now we show that $\mathcal{P}_{\mathcal{F}}= L^{*}_{t,s}:f$. Suppose that $\mathcal{P}_{\mathcal{F}}$ is a proper subset of $L^{*}_{t,s}:f$.
Hence there exists a monomial $f'$, in the variables $y_{r}$ with $r\notin \mathcal{F}$, of degree at least 1 such that $ff'\in L^{*}_{t,s}$.
This implies that there exists a monomial $u=y_{1}^{b_{1}}\dots y_{m_{1}+\dots+m_{n}}^{b_{m_{1}+\dots+m_{n}}}\in G(L^{*}_{t,s})$ such that $u$ divides
$ff'$. Thus $b_{r}\leq d_{r}$ for any $r\in \mathcal{F}$ because $f'\in K[y_{r}\mid r\notin \mathcal{F}]$.
It follows that
\[
t=\deg(u)=\Sigma_{r=1}^{m_{1}+\dots+m_{n}}b_{r}\leq \Sigma_{r\in \mathcal{F}}d_{r}+s(m_{1}+\dots+m_{n}-|\mathcal{F}|)=\deg(f)=t-1,
\]
which is a contradiction. Therefore, $\mathcal{P}_{\mathcal{F}}$ is not a proper subset of $L^{*}_{t,s}:f$, but $\mathcal{P}_{\mathcal{F}}=L^{*}_{t,s}:f$.
This equality means that $\mathcal{P}_{\mathcal{F}}\in \Ass_{T}(T/L^{*}_{t,s})$.
\end{proof}

\begin{example}
\label{L_{15,4}}
{\em
Let
\begin{eqnarray*}
L^{*}_{15,4}=(x_{11}^{4}x_{12}^{4}x_{21}^{4}x_{22}^{3}, x_{11}^{4}x_{12}^{4}x_{21}^{3}x_{22}^{4},
x_{11}^{4}x_{12}^{3}x_{21}^{4}x_{22}^{4}, x_{11}^{3}x_{12}^{4}x_{21}^{4}x_{22}^{4}),
\end{eqnarray*}
be a generalized Veronese bi-type ideal of $T=K[x_{11},x_{12},x_{21},x_{22}]$.
Therefore, Theorem ~\ref{Ass} yields
\begin{eqnarray*}
\Ass_{T}(T/L^{*}_{15,4})&=&\{(x_{11}), (x_{12}), (x_{21}), (x_{22}), (x_{11},x_{12}), (x_{11},x_{21}), (x_{11},x_{22}),
(x_{12}, x_{21}), \\
&&(x_{12},x_{22}), (x_{21}, x_{22})\}.
\end{eqnarray*}
}
\end{example}

\section{Toric ideal of $K[L^{*}_{t,s}]$}
\label{four}
Let $T=K[x_{11},\ldots,x_{1m_1},x_{21},\ldots,x_{2m_2},\ldots,x_{n1},\ldots,x_{nm_n}]$ be a polynomial ring over a field $K$ in the variables
\[
x_{11},\ldots,x_{1m_1},x_{21},\ldots,x_{2m_2},\ldots,x_{n1},\ldots,x_{nm_n},
\]
and let $L^{*}_{t,s}=(f_{1},\dots,f_{p})$ be the ideal of generalized Veronese bi-type. The monomial subring of $T$ spanned by $F=\{f_{1},\dots,f_{p}\}$
is the $K$-algebra $K[L^{*}_{t,s}]=K[F]=K[f_{1},\dots,f_{p}]$.
The monomial subring $K[F]$ is a graded subring of $T$ with the grading given by $K[F]_{i}=K[F]\sect T_{i}$.
There is a graded epimorphism of $K$-algebras:
\[
\varphi: R=K[t_{1},\dots, t_{p}]\rightarrow K[L^{*}_{t,s}] \rightarrow 0, \quad \text{induced by} \quad \varphi(t_{i})=f_{i},
\]
where $R$ is a graded polynomial ring with the grading induced by setting $\deg(t_{i})=\deg(f_{i})$.
Notice that the map $\varphi$ is given by $\varphi(h(t_{1},\dots,t_{p}))=h(f_{1},\dots,f_{p})$ for all $h\in R$.

The kernel of $\varphi$, denoted by $P_{F}$, is the so-called {\em toric ideal} of $K[L^{*}_{t,s}]$ with respect to $f_{1},\dots,f_{p}$.
We also denote the toric ideal of $K[L^{*}_{t,s}]$ by $I(L^{*}_{t,s})$.

In this section we prove that $I(L^{*}_{t,s})$ has a quadratic Groebner basis. In order to formulate this result we have to recall
the notion sortability, introduced \cite{S}.

Let $A=K[z_{1},\dots, z_{q}]$ be a polynomial ring and $L$ be a monomial ideal of $A$ generated in degree $t$.
Let $\mathcal{B}$ be the set of the exponent vectors of the monomials of $G(L)$. If $u=(u_{1},\dots,u_{q})$, $v=(v_{1},\dots,v_{q})\in \mathcal{B}$,
then $\underline{z}^{u}=\prod_{i}z_{i}^{u_{i}}$, $\underline{z}^{v}=\prod_{i}z_{i}^{v_{i}}\in L$.
Then we write $\underline{z}^{u}\underline{z}^{v}=z_{i_{1}}\dots z_{i_{2t}}$ with $i_{1} \leq i_{2} \leq \dots \leq i_{2t}$. We set
$\underline{z}^{u'}=\prod _{j=1}^{t}z_{2j-1}$ and $\underline{z}^{v'}=\prod _{j=1}^{t}z_{2t}$.
This defines a map
\[
sort: \mathcal{B}\times\mathcal{B}\rightarrow M_{t}\times M_{t}, \quad (u,v)\rightarrow (u',v'),
\]
where $M_{t}$ is the set of all integer vectors $(a_{1},\dots,a_{q})$ such that $\sum_{i=1}^{q}a_{i}=t$.
The set $\mathcal{B}$ is called {\em sortable} if $\Im(sort)\subseteq \mathcal{B}\times\mathcal{B}$.

The ideal $L$ is called {\em sortable} if the set of exponent vectors of the monomials of $G(L)$ is sortable.
In other words, let $\underline{z}^{u}$, $\underline{z}^{v}\in L$, then $L$ is said sortable if $\underline{z}^{u'}$, $\underline{z}^{v'}\in L$, where
$(u',v')=sort(u,v)$.

\begin{theorem}
\label{sort}
Let $L^{*}_{t,s}$ be a generalized Veronese bi-type ideal of $T$. Then $L^{*}_{t,s}$ is sortable.
\end{theorem}

\begin{proof}
Let $T=K[x_{11},\ldots,x_{1m_1},x_{21},\ldots,x_{2m_2},\ldots,x_{n1},\ldots,x_{nm_n}]$ be a polynomial ring over a field $K$ in the variables
$x_{11},\ldots,x_{1m_1},x_{21},\ldots,x_{2m_2},\ldots,x_{n1},\ldots,x_{nm_n}$,
and let
\[
L^{*}_{t,s}=(\{x_{11}^{a_{11}} \dots x_{1m_1}^{a_{1m_{1}}} \dots x_{n1}^{a_{n1}} \dots x_{nm_n}^{a_{nm_{n}}}| \sum_{i=1}^{n}\sum_{j=1}^{m_{i}}a_{ij}=t, \quad  0\leq a_{ij}\leq s\}),
\]
be the ideal of generalized Veronese bi-type.
Furthermore, let $\mathcal{B}$ be the set of the exponent vectors of the monomials of $G(L^{*}_{t,s})$.
Let $f_{i}=x_{11}^{a_{11}} \dots x_{1m_1}^{a_{1m_{1}}} \dots x_{n1}^{a_{n1}} \dots x_{nm_n}^{a_{nm_{n}}}$, $f_{j}=x_{11}^{b_{11}} \dots x_{1m_1}^{b_{1m_{1}}} \dots x_{n1}^{b_{n1}} \dots x_{nm_n}^{b_{nm_{n}}}\in G(L^{*}_{t,s})$, then
\[
u=(a_{11},\ldots,a_{1m_1};\ldots;a_{n1},\ldots,a_{nm_n}),
v=(b_{11},\ldots,b_{1m_1};\ldots;b_{n1},\ldots,b_{nm_n})\in \mathcal{B}.
\]
Therefore we obtain that
\[
f_{i}f_{j}=\underbrace{x_{11}\dots x_{11}}_{a_{11}+b_{11}-times}\dots \underbrace{x_{1m_1}\dots x_{1m_1}}_{a_{1m_1}+b_{1m_1}-times} \dots \underbrace{ x_{n1}\dots x_{n1}}_{a_{n1}+b_{n1}-times}\dots \underbrace{x_{nm_n}\dots x_{nm_n}}_{a_{nm_n}+b_{nm_n}-times}
\]
is a monomial of degree $2t$. If one replaces the set of variables $\{x_{11},\dots, x_{1m_{1}}\}$ with $\{z_{1},\dots,z_{m_{1}}\}$ and $\{x_{21},\dots,x_{2m_{2}}\}$ with $\{z_{m_{1}+1},\dots, z_{m_{1}+m_{2}}\}$ and so on up to $\{x_{n1},\dots,x_{nm_{n}}\}$ with
$\{z_{m_{1}+\dots +m_{n-1}+1},\dots,z_{m_{1}+\dots+m_{n}}\}$, thus $f_{i}f_{j}=z_{i_{1}}\dots z_{i_{2t}}$ with $i_{1}\leq \dots \leq i_{2t}$.
We consider $f_{i}'=\underline{z}^{u'}=\prod _{r=1}^{t}z_{2r-1}$ and $f_{j}'=\underline{z}^{v'}=\prod_{r=1}^{t}z_{2t}$.
We prove that $f_{i}',f_{j}'\in L^{*}_{t,s}$. Observe that $f_{i}'$ is of degree $t$ and we write
\[
f_{i}'=\prod _{r=1}^{t}z_{2r-1}=x_{11}^{a'_{11}} \dots x_{1m_1}^{a'_{1m_{1}}} \dots x_{n1}^{a'_{n1}} \dots x_{nm_n}^{a'_{nm_{n}}}.
\]
If $a_{ij}+b_{ij}$ is even then $a_{ij}'=\frac{a_{ij}+b_{ij}}{2}\leq s$ and if  $a_{ij}+b_{ij}$ is odd then $a_{ij}'=\frac{a_{ij}+b_{ij}+1}{2}< s$.
Furthermore, because $f_{i}'$ is of degree $t$ and there exists $a_{ij}'\neq 0$ with $0\leq a_{ij}' \leq t$ for all $ij$, thus
$x_{i1}^{a'_{i1}} \dots x_{im_1}^{a'_{im_{i}}}\neq x_{ij}^{t}$. This implies that $x_{i1}^{a'_{i1}} \dots x_{im_1}^{a'_{im_{i}}}\in L_{i,q_{i},s}$ for all $i=1,\dots, n$ with
$\sum_{i=1}^{n}q_{i}=t$. Therefore, $f_{i}'\in L^{*}_{t,s}$. In the same way the argument holds for $f_{j}'$, and hence $L^{*}_{t,s}$ is sortable.
\end{proof}

\begin{corollary}
\label{Groebner}
Let $L^{*}_{t,s}$ be a generalized Veronese bi-type ideal of $T$. Then:

\item[(1)] $I(L^{*}_{t,s})$ has a quadratic Groebner basis.

\item[(2)] $K[L^{*}_{t,s}]$ is Koszul.
\end{corollary}

\begin{proof}
(1) The assertion follows by Theorem \ref{sort} and \cite[Lemma 5.2]{HH}.

(2) The conclusion follows by (1).
\end{proof}

%%}
\label{'ubl'}
\end{document}